\documentclass[preprint]{elsarticle}
\usepackage[margin=24mm]{geometry}
\usepackage{amsmath,amssymb,amsthm,bm,graphicx,booktabs,array}
\usepackage[T1]{fontenc}
\usepackage{lmodern,microtype}
\usepackage{float}
\usepackage[colorlinks=true,linkcolor=blue,citecolor=blue,urlcolor=blue]{hyperref}
\makeatletter
\long\def\@makecaption#1#2{\vskip\abovecaptionskip
  {\small\noindent\textbf{#1:} #2\par}\vskip\belowcaptionskip}
\makeatother
\let\originalthebibliography\thebibliography
\renewcommand{\thebibliography}[1]{\footnotesize\originalthebibliography{#1}}
\newcommand{\C}{\mathbb C}
\newcommand{\R}{\mathbb R}
\newcommand{\inc}{\mathrm{inc}}
\newcommand{\exc}{\mathrm{exc}}
\newcommand{\GT}{\mathrm{ref}}
\newcommand{\out}{\mathrm{out}}
\newcommand{\norm}[1]{\left\lVert#1\right\rVert}

\journal{Computers Physics Communications}

\begin{document}
\begin{frontmatter}

\title{ESPC-Transolver: An Equivalent-Surface Physics-Guided Coupling-Aware Transolver for 3D Multiple Scattering}

\author[lab1]{Guoqing Zhao}
\ead{lqing2781@gmail.com}
\author[lab1]{Liang Li\corref{cor1}} 
\ead{plum\_liliang@uestc.edu.cn, plum.liliang@gmail.com}
\author[lab2]{Xing Li\corref{cor1}} 
\ead{xingli@uestc.edu.cn}
\author[lab1]{Hao Kang}
\ead{kanghao2022@gmail.com}

\address[lab1]{School of Mathematical Sciences, University of
               Electronic Science and Technology of China, 
               Chengdu, P.R. China}

\address[lab2]{School of Electronic Science and Engineering , University of
               Electronic Science and Technology of China, 
               Chengdu, P.R. China}

\cortext[cor1]{Corresponding authors}

\begin{abstract}
Cross-geometry generalization is particularly challenging because geometry-conditioned multiple-scattering surrogates must infer unseen objects' scattering properties from geometric features while accurately capturing complex, phase-sensitive interactions arising from repeated wave exchange. We introduce ESPC-Transolver\footnote{Code repository: \url{https://github.com/xiaozhao111111111/ESPC-Transolver-code}}, an equivalent-surface physics-guided, coupling-aware Transolver that separates local scattering-response representations from inter-object interaction modeling. For each sound-hard scatterer, canonical scattered-pressure responses are constructed under a shared set of spherical-wave excitations. A local Dirichlet-to-Neumann (DtN) map and modal scattering conversion express these responses on an enclosing sphere, providing physically consistent input channels across geometries. Shared object-local Transolver blocks encode these responses and incident-field traces, while relative-position-aware slice attention models inter-object coupling. Predicted interface pressure and normal flux recovered from the local operator define equivalent sources for reconstructing exterior near- and far-field scattering. We evaluate the method on three-dimensional acoustic multiple-scattering problems at a fixed frequency, comparing it with geometry-conditioned predictors and alternative neural network backbones. ESPC-Transolver achieves lower scattered-field prediction errors on held-out instances of familiar geometry families and generalizes effectively to previously unseen Peanut and Torus geometries. Exterior-field reconstructions reproduce multiple-scattering interference patterns, linking cross-geometry interface prediction to near-field pressure and far-field scattering recovery.
\end{abstract}
\begin{keyword} 
Canonical boundary response\sep
Cross-geometry generalization\sep
Equivalent sources \sep
Transolver
\end{keyword}
\end{frontmatter}

\section{Introduction}
Multiple scattering occurs when the field scattered by one object modifies the excitation incident on other objects, giving rise to a coupled wave field through repeated interactions. The resulting field depends not only on the geometry of individual scatterers but also on their relative positions, orientations, and the phase relationships between interacting waves \cite{martin2006}. Boundary element methods (BEM) formulate exterior scattering problems through boundary integral equations defined on the surfaces of the scatterers \cite{bonnet1999}. By employing the outgoing Green's function, these methods satisfy the radiation condition without discretizing the unbounded exterior domain. The resulting coupled boundary system accounts for multiple interactions among scatterers, while fast multipole methods reduce the computational cost of evaluating dense interactions \cite{darve2000}. Nevertheless, changes in the geometry, arrangement, or incident excitation generally require the corresponding coupled system to be solved again. Deep-learning surrogates that predict the boundary data needed for field reconstruction offer an alternative for repeated scattering calculations, while retaining the physical relationship between boundary quantities and the exterior field.

Deep-learning approaches to wave scattering have been developed at different levels of representation. Direct field surrogates learn mappings from scattering configurations to measured or spatially sampled fields. SwitchNet exploits the low-rank structure of scattering interactions to construct sparse switching layers for forward and inverse wave-scattering mappings \cite{khoo2019}. Fourier neural operators (FNO) and DeepONet provide general frameworks for learning operators that map partial differential equation inputs to their corresponding solutions \cite{li2020,lu2021}. For oscillatory wave fields, MscaleFNO combines Fourier operators acting on differently scaled inputs to improve the representation of high-frequency components in Helmholtz coefficient-to-solution mappings \cite{you2026}. SwitchNet and MscaleFNO thus illustrate how scattering-specific architectures exploit interaction structure and oscillatory scales, complementing the general operator-learning frameworks.

Geometry-conditioned operator learning extends this setting to problems in which the scattering boundary varies. Nair et al. combined a physics-informed DeepONet with a non-uniform rational B-spline (NURBS) parameterization to predict acoustic scattered pressure for rigid obstacles with varying geometries \cite{nair2025}. Their approach illustrates how explicit geometric representations can support learning across changing computational domains. In multiple scattering, however, the predictor must also account for the phase and amplitude of waves exchanged among objects. The choice of geometric conditioning variables therefore affects both the representation of individual scattering responses and the prediction of scene-dependent interactions.

Boundary-based learning provides a complementary approach to scattering surrogates. The boundary integral network with an unknown Green's function (BIN-G) learns a domain-independent Green's function and associated boundary density functions using training conditions derived from partial differential equations and boundary integral formulations, with numerical demonstrations on variable-coefficient Poisson and Helmholtz problems \cite{negi2024}. PIBNet uses a multiscale graph neural network to predict boundary unknowns and represents interactions between distant obstacles through graph edges \cite{marsal2025}. Neural operators can also be integrated into conventional numerical solvers rather than used to predict the complete field directly. Lee et al. employed a DeepONet trained on a domain without scatterers to construct coarse-scale preconditioners for iterative three-dimensional scattering solvers \cite{lee2025}. These approaches illustrate distinct roles for learning in scattering computations, including the approximation of propagation operators, prediction of boundary solutions, and acceleration of coupled linear systems.

Attention-based architectures offer a flexible mechanism for modeling interactions across irregular spatial discretizations. Transolver aggregates spatial points into learned slices and applies Physics-Attention to model interactions between the resulting representations \cite{wu2024}. Transolver++ extends this framework to larger geometric configurations \cite{luo2025}. These architectures provide a general mechanism for exchanging information across spatial locations. For multiple scattering, the effectiveness of such interactions also depends on the physical information represented by the input features. Coordinates and surface normals describe the geometry of a scattering boundary, whereas a local scattering-response operator characterizes how the object transforms an incident excitation into a scattered field. This distinction motivates the use of scattering responses, rather than geometry alone, as the information exchanged within a learning-based multiple-scattering model.

These approaches differ in what they learn: shape-conditioned fields, propagation kernels, boundary unknowns, or corrections within a numerical solver. For a geometry family absent from training, a geometry-conditioned predictor must infer the object's scattering response and incorporate it into phase-sensitive interactions with other objects. Responses to a common set of wave excitations instead give each input channel a shared physical meaning, allowing object-specific scattering information to enter the coupling model directly. This raises a specific question: does conditioning a coupling-aware network on such local responses improve prediction across unseen geometry families compared with geometric conditioning?

We introduce ESPC-Transolver, an equivalent-surface physics-guided, coupling-aware Transolver for three-dimensional acoustic multiple scattering. The method represents individual scatterers through canonical scattering responses and uses a shared local encoder together with coupling-aware attention to predict their interactions in multi-object configurations. Building on local scattering operators and artificial-interface formulations \cite{martin2006,pedneault2017}, the method makes three contributions:

\begin{enumerate}
\setlength{\itemsep}{3pt}
\setlength{\parskip}{0pt}
\item \textbf{Canonical scattering-response representation for cross-geometry learning.} We construct canonical scattered-pressure responses for sound-hard objects using a common set of regular spherical-wave excitations. A local Dirichlet-to-Neumann (DtN) map and modal scattering conversion express these responses on enclosing spheres. By assigning consistent physical meanings to the excitation modes and observation channels, the representation describes geometry-dependent scattering properties within a common response space. It thereby provides physically aligned inputs for learning across different geometry families.

\item \textbf{A coupling-aware Transolver architecture.} We develop shared object-local Transolver blocks that encode canonical scattering responses and incident-field traces. Relative-position-aware slice attention, incorporating spatial and wave-dependent features, enables information exchange between objects and predicts the scattered-pressure traces under multiple-scattering interactions. The architecture distinguishes local response encoding from inter-object coupling while sharing parameters across scatterers, allowing configurations with different numbers and arrangements of objects to be processed within a unified framework.

\item \textbf{Equivalent-interface learning and exterior-field reconstruction across geometries and topologies.} We represent scatterers with different shapes and surface topologies on common enclosing spherical interfaces. Geometric and topological variations are incorporated through the local scattering-response operators, while the canonical sampling layout and coupling network remain unchanged. The predicted interface pressure, together with the normal flux obtained from the local operator, defines equivalent sources for reconstructing the exterior near-field pressure and far-field scattering response through a common procedure.
\end{enumerate}

We evaluate the method on three-dimensional acoustic multiple-scattering problems at a fixed frequency. First, canonical scattering responses are compared with radial-occupancy features and surface point-and-normal representations under a shared network backbone to assess the effect of input representation. We then compare alternative interaction architectures using the same canonical-response inputs to examine the role of coupling-aware network design. Generalization is evaluated on held-out instances of familiar geometry families and on previously unseen Peanut and Torus families. Representative near-field and far-field reconstructions further assess how errors in the predicted interface quantities affect the recovered spatial distributions and interference patterns of the scattered field.

The remainder of this paper is organized as follows. Section~\ref{sec:response} develops the canonical scattering-response formulation. Section~\ref{sec:model} presents the coupling-aware Transolver architecture. Section~\ref{sec:recovery} describes equivalent-source reconstruction of the exterior scattered field. Section~\ref{sec:results} reports the numerical experiments and comparative results.
\section{Problem formulation and canonical equivalent interfaces}
\label{sec:response}

\subsection{The exterior scattering problem}

\begingroup
\predisplaypenalty=3000
Let \(\Omega_i\subset\R^3\), \(i=1,\ldots,N\), be sound-hard objects centred at \(c_i\), and let \(D^e=\R^3\setminus\bigcup_i\overline{\Omega_i}\). The total pressure \(u=u^\inc+u^s\) satisfies
\begin{align}
 \Delta u+k^2u&=0 &&\text{in }D^e,\label{eq:pde}\\
 \partial_nu&=0 &&\text{on }\partial\Omega_i,\label{eq:neumann}\\
 u^\inc(x)&=\exp(ikd\cdot x),\qquad |d|=1,\label{eq:inc}\\
 \lim_{r\to\infty}r(\partial_ru^s-iku^s)&=0,\label{eq:radiation}
\end{align}
\endgroup
where \(k\) is the wavenumber, \(d\) is the incident direction and \(u^s\) is the outgoing scattered pressure. The outward normal \(n\) on an obstacle defines the sound-hard condition. Spherical Hankel functions \(h_l^{(1)}\) follow the outgoing convention in Eq.~\eqref{eq:radiation}.

\begin{figure}[htbp]
\centering
\includegraphics[width=0.70\textwidth]{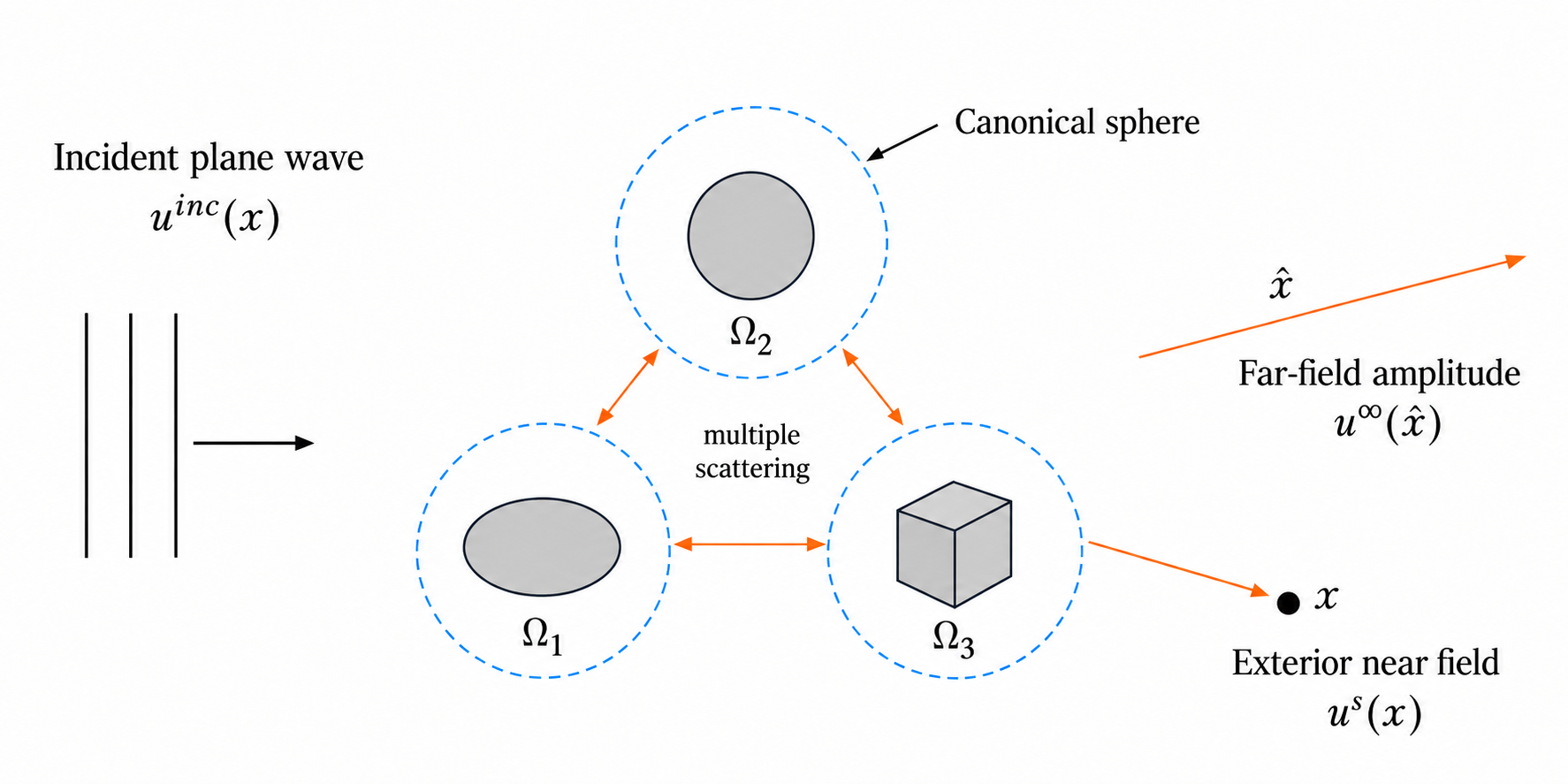}
\caption{Sound-hard multiple scattering and exterior field recovery. An incident plane wave illuminates three objects. Blue dashed outlines mark disjoint canonical spheres, and bidirectional arrows indicate wave exchange between objects. Canonical pressure and flux traces define equivalent sources for the scattered near field outside all canonical balls and the far-field amplitude. Each recovered field includes contributions from every object.}
\label{fig:scattering_problem}
\end{figure}

Each object is enclosed by \(\Gamma_i^c=\partial B_{R_c}(c_i)\), with the same physical radius \(R_c\). The canonical balls are disjoint and each contains its designated object alone. Translation to object-centred coordinates changes neither the object's scale nor the physical wavenumber. Figure~\ref{fig:scattering_problem} shows the physical objects, their canonical interfaces and the exterior observation regions.

\subsection{A local response on the canonical interface}

Within \(D_i^c=B_{R_c}(c_i)\setminus\overline{\Omega_i}\), prescribe \(u=g\) on \(\Gamma_i^c\) and impose the sound-hard condition \(\partial_nu=0\) on \(\partial\Omega_i\). Assume a Lipschitz obstacle boundary strictly enclosed by the canonical sphere. Let \(\mathcal A_i=-\Delta\) denote the mixed Laplacian with homogeneous Neumann data on the obstacle and homogeneous Dirichlet data on the canonical sphere. For \(k^2\notin\sigma(\mathcal A_i)\), the local problem is uniquely solvable and defines a bounded DtN map
\begin{equation}
 \Lambda_i:g\longmapsto q=\left.\partial_{n_c}u\right|_{\Gamma_i^c},
 \label{eq:dtn}
\end{equation}
where \(n_c\) points radially outward from the canonical ball. The operator depends on the object's geometry and physical wavenumber, but not on its scene position or incident direction. The restriction of a coupled scene field to this local domain satisfies \(q_i=\Lambda_i u_i\), transferring the object's boundary response to a common interface.

Choose spherical-harmonic probes
\begin{equation}
 V_Q=\operatorname{span}\{Y_l^m:0\leq l\leq L,\ -l\leq m\leq l\},
 \qquad Q=(L+1)^2.
 \label{eq:basis}
\end{equation}
At \(M\) common directions, let \(\Phi_Q\in\C^{M\times Q}\) contain the sampled harmonics and let \(W\) be the positive diagonal matrix of sphere-area quadrature weights. Full column rank of \(\Phi_Q\) makes its weighted Gram matrix invertible. The weighted projection of a sampled trace to modal coefficients is
\begin{equation}
 P_Q=(\Phi_Q^*W\Phi_Q)^{-1}\Phi_Q^*W.
 \label{eq:projection}
\end{equation}
It satisfies \(P_Q\Phi_Q=I_Q\) on the sampled harmonic space. The probes define a common input basis and the canonical directions define a common output grid. A finite element method (FEM) constructs the sampled response by projecting weak reactions onto output harmonics and evaluating the resulting coefficients on this grid.

\subsection{Constructing the discrete response by internal elimination}

A curved quadratic FEM mesh discretizes \(D_i^c\). Let \(V_g=\{u\in H^1(D_i^c):u|_{\Gamma_i^c}=g\}\) and \(V_0=\{v\in H^1(D_i^c):v|_{\Gamma_i^c}=0\}\). With the sesquilinear form linear in its first argument, the weak problem is~\cite{brenner2008}
\begin{equation}
 a(u,v)=\int_{D_i^c}\left(\nabla u\cdot\nabla\overline v-k^2u\overline v\right)\,dx=0,\qquad v\in V_0.
 \label{eq:weakform}
\end{equation}
The homogeneous Neumann condition on the obstacle enters naturally. For a test trace \(w\in H^{1/2}(\Gamma_i^c)\), choose any FEM lifting \(v_w\) with trace \(w\) on \(\Gamma_i^c\). The weak conormal reaction is defined by \(\langle q,w\rangle_{\Gamma_i^c}=a(u,v_w)\), with \(\langle q,w\rangle=\int_{\Gamma_i^c}q\,\overline w\,dS\) for smooth traces. Let \(b\) index the actual FEM degrees of freedom on the canonical boundary and \(I\) all remaining degrees of freedom, including those on the Neumann obstacle boundary. Partition the discrete Helmholtz system as
\begin{equation}
 \begin{bmatrix}A_{bb}&A_{bI}\\A_{Ib}&A_{II}\end{bmatrix}
 \begin{bmatrix}U_b\\U_I\end{bmatrix}
 =\begin{bmatrix}R_b\\0\end{bmatrix}.
 \label{eq:femblocks}
\end{equation}
The vector \(R_b\) stores weak reactions paired with the FEM trace basis. Eliminating \(U_I\) gives
\begin{equation}
 U_I=-A_{II}^{-1}A_{Ib}U_b,\qquad
 R_b=\left(A_{bb}-A_{bI}A_{II}^{-1}A_{Ib}\right)U_b.
 \label{eq:schur}
\end{equation}
Thus the Schur action is the discrete DtN response on the FEM trace space~\cite{brenner2008}. One factorization of \(A_{II}\) is reused for the probe solves.

Let \(B_{\mathrm{in}}\in\C^{n_b\times Q_{\mathrm{in}}}\) contain the FEM trace coefficients of the spherical-harmonic probes. The consistent boundary mass matrix \(M_b\in\C^{n_b\times n_b}\) maps nodal flux coefficients to weak reactions. The reaction projection is the weighted least-squares fit
\(C_i=\arg\min_C\norm{q_b-H_{\mathrm{out}}C}_{M_b}^2\), whose normal equations are shown below. The multiple-right-hand-side solve and reaction projection are
\begin{align}
 X&=-A_{II}^{-1}A_{Ib}B_{\mathrm{in}},\qquad
 R=A_{bb}B_{\mathrm{in}}+A_{bI}X,\qquad M_bq_b=R,\\
 G_b&=H_{\mathrm{out}}^*M_bH_{\mathrm{out}},\qquad
 G_bC_i=H_{\mathrm{out}}^*R,\qquad
 Z_i=\Phi_{\mathrm{out}}C_iE_Q.
 \label{eq:fluxprojection}
\end{align}
Here \(H_{\mathrm{out}}\in\C^{n_b\times Q_{\mathrm{out}}}\) contains the output harmonic traces, \(G_b\) is positive definite when these traces are linearly independent in the \(M_b\)-weighted inner product, and \(C_i\in\C^{Q_{\mathrm{out}}\times Q_{\mathrm{in}}}\) contains the projected flux coefficients. The selection matrix \(E_Q\in\R^{Q_{\mathrm{in}}\times Q}\) retains the input probes used for conditioning, and \(\Phi_{\mathrm{out}}\in\C^{M\times Q_{\mathrm{out}}}\) evaluates the projected flux on the common canonical grid. Thus the construction proceeds from FEM weak reactions to mass-matrix modal projection and finally to canonical sampling. Appendix~\ref{app:numerics} specifies the input, output and sampling dimensions.

\subsection{Converting the DtN response to scattered pressure}

Converting the DtN response to scattered pressure aligns the conditioning representation with the prediction target. Let \(B_i=P_QZ_i\in\C^{Q\times Q}\). Write the retained canonical traces in terms of regular incoming coefficients \(a_i^{\mathrm{exc}}\) and outgoing coefficients \(b_i\)~\cite[Chapter~7]{martin2006}:
\begin{equation}
 u_i=\Phi_Q(Ja_i^{\mathrm{exc}}+Hb_i),\qquad
 q_i=\Phi_Q(J'a_i^{\mathrm{exc}}+H'b_i).
 \label{eq:localwaves}
\end{equation}
The diagonal matrices \(J\) and \(H\) in \(\C^{Q\times Q}\) contain \(j_l(kR_c)\) and \(h_l^{(1)}(kR_c)\), repeated over \(m\). Their primes denote derivatives with respect to physical radius, including the factor \(k\). Applying the projected DtN relation to Eq.~\eqref{eq:localwaves} gives
\begin{equation}
 B_i(Ja_i^{\mathrm{exc}}+Hb_i)=J'a_i^{\mathrm{exc}}+H'b_i,\qquad
 (H'-B_iH)b_i=(B_iJ-J')a_i^{\mathrm{exc}}.
 \label{eq:conversion}
\end{equation}
When \(H'-B_iH\) is nonsingular, the outgoing response operator and canonical scattering response are
\begin{equation}
 T_i^\out=(H'-B_iH)^{-1}(B_iJ-J'),\qquad
 S_i(Z)=\Phi_QHT_i^\out.
 \label{eq:sresponse}
\end{equation}
Here \(T_i^{\out}\in\C^{Q\times Q}\) and \(S_i(Z)\in\C^{M\times Q}\). Each column of \(S_i(Z)\) is the canonical scattered-pressure trace produced by one regular spherical-wave excitation. The response is defined by the FEM discretization, modal spaces and quadrature, with \(P_Q\Phi_Q=I_Q\) on the represented sampled space. ESPC-Transolver receives this response together with the incident trace and scene features.

For the retained modal system, let \(a_i^{\inc}\) denote the plane-wave coefficients expanded in regular spherical waves about centre \(c_i\). The total regular excitation is \(a_i^{\exc}=a_i^{\inc}+\sum_{j\ne i}U_{ij}b_j\), where the off-diagonal translation blocks map outgoing waves at one centre to regular waves at another~\cite[Chapter~7]{martin2006}. With \(\mathsf T=\operatorname{diag}(T_1^{\out},\ldots,T_N^{\out})\), the Foldy--Lax system is
\begin{equation}
 b=\mathsf T(a^{\inc}+\mathsf U b),\qquad (I-\mathsf T\mathsf U)b=\mathsf T a^{\inc}.
 \label{eq:foldylax}
\end{equation}
The diagonal blocks of \(\mathsf U\) are zero, while each object's self response is contained in \(T_i^{\out}\). The retained coupled system is well-defined when \(I-\mathsf T\mathsf U\) is nonsingular. ESPC-Transolver learns the scene-to-trace map
\begin{equation}
 \mathcal N_\theta\left(\{S_i,\mathbf c_i,u_i^{\inc}\}_{i=1}^{N}\right)
 =\{\widehat u_i^s\}_{i=1}^{N},
 \label{eq:networkmap}
\end{equation}
using local canonical responses and scene context to predict the coupled scattered-pressure traces. The Foldy--Lax equation describes the retained modal coupling, while the network predicts its canonical pressure-trace output. Appendix~\ref{app:theory} derives the residual relation between the recovered outgoing coefficients and this coupling equation, together with bounds on the propagation of local-response perturbations through the coupled system.

\section{ESPC-Transolver: a compositional Transolver architecture}
\label{sec:model}

ESPC-Transolver separates processing on each canonical sphere from communication between objects (Fig.~\ref{fig:method}). Each object is represented by its canonical scattering response, incident trace, local coordinates and illumination features. Shared local blocks encode these inputs, while an interaction module uses relative positions to exchange information between objects.

\subsection{Object-compositional organization}

\begin{figure}[!htb]
\centering
\includegraphics[trim=0 0 0 79bp,clip,width=\textwidth]{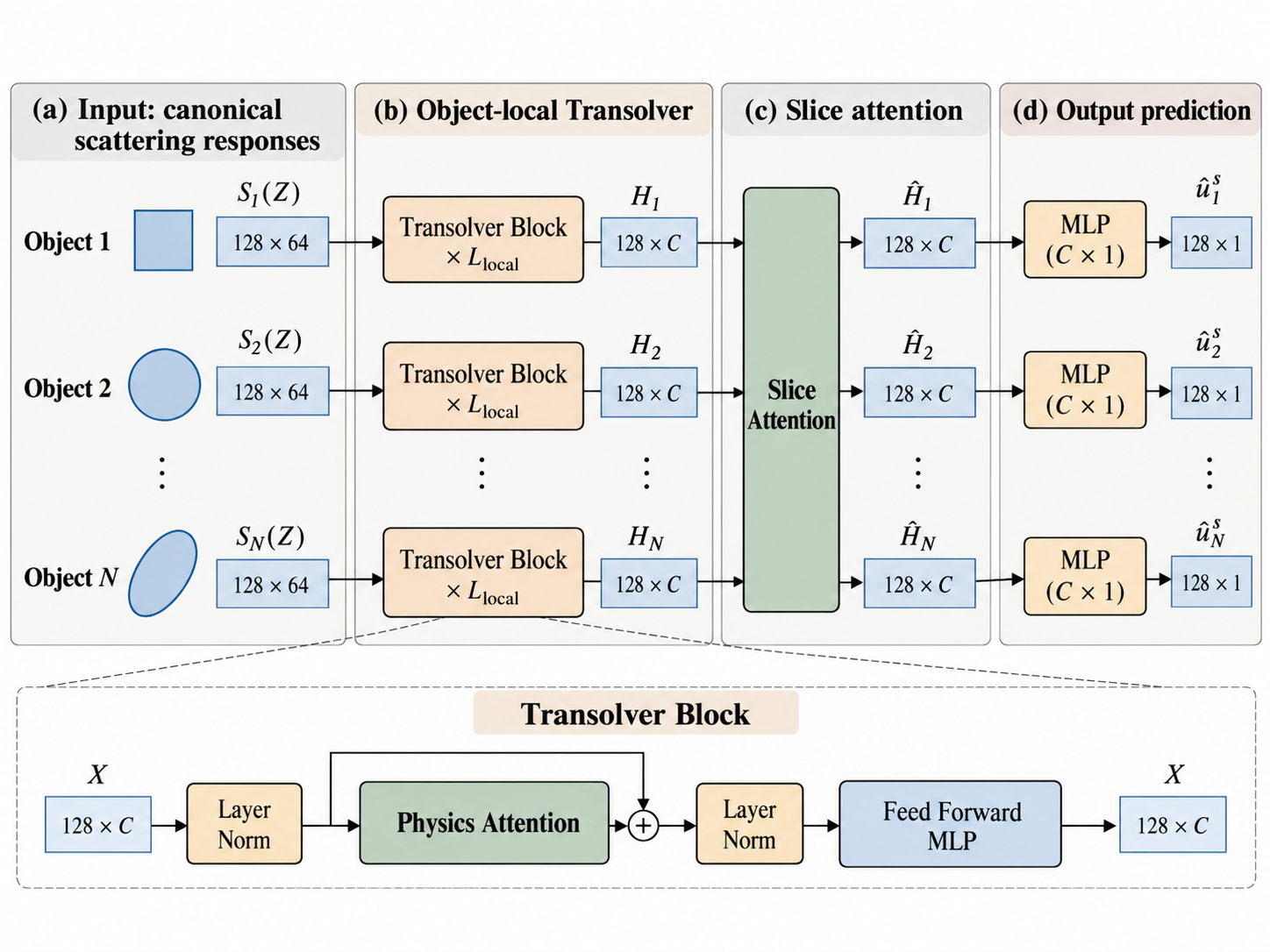}
\caption{ESPC-Transolver architecture. (a) Canonical scattering responses condition each object. (b) Shared object-local Transolver blocks process its canonical nodes. (c) Relative-position-aware slice attention exchanges information between objects. Local and inter-object updates are interleaved through the network depth. (d) A shared decoder produces one complex scattered-pressure value per node, represented by two real channels. The lower inset shows the Physics-Attention and feed-forward structure of a local block. Incident traces and illumination features accompany the response input.}
\label{fig:method}
\end{figure}

For object \(i\), the input combines the real and imaginary parts of \(S_i(Z)\) with the incident trace and local node features. A shared embedding maps them to \(H_i\in\R^{M\times C}\), where \(M\) is the number of canonical nodes and \(C\) is the hidden width. Object-centred coordinates enter the local representation; object centres are used only to construct relative positions for inter-object exchange. Figure~\ref{fig:method} summarizes this data flow.

Each object is processed independently by shared Transolver blocks using Physics-Attention \cite{wu2024}. For head \(h\), learned slice logits define assignments \(\pi_{irs}^{(h)}=\operatorname{softmax}_{s}(\ell_{ir}^{(h)})\), normalized over slices for each node \(r\). These assignments aggregate nodes into slice tokens, attention updates the tokens, and the same assignments return the updates to the canonical nodes. Layer normalization, attention and feed-forward sublayers use residual updates. The local-block parameters are shared across objects, while successive blocks have distinct parameters.

The local backbone has five blocks, each with eight heads of width \(d_h=32\) and 32 slices per head. For a local head, node values \(v_{ir}^{(h)}\) and assignments \(\pi_{irs}^{(h)}\) form normalized slice tokens \(t_{is}^{(h)}\in\R^{d_h}\). Attention over those tokens gives
\begin{align}
 t_{is}^{(h)}&=\frac{\sum_r\pi_{irs}^{(h)}v_{ir}^{(h)}}{\sum_r\pi_{irs}^{(h)}},&
 \widetilde t_{is}^{(h)}&=\sum_t
 \operatorname{softmax}_t\!\left(\frac{(q_{is}^{(h)})^{\top}k_{it}^{(h)}}{\sqrt{d_h}}\right)v_{it}^{\mathrm{slice},(h)},
 \label{eq:localattention}\\
 \widetilde F_{ir}&=W_o\operatorname{concat}_{h}\left(\sum_s\pi_{irs}^{(h)}\widetilde t_{is}^{(h)}\right).
 \label{eq:localdeslice}
\end{align}
Here \(v_{ir}^{(h)}\) is a projected node value, while \(q_{is}^{(h)}\), \(k_{is}^{(h)}\) and \(v_{is}^{\mathrm{slice},(h)}\) are token queries, keys and values. The softmax is over key slices for each target slice. Heads are concatenated and projected to \(C=256\), and the block's feed-forward residual completes the within-object update.

An interaction module is inserted after local blocks 1, 2 and 3. It forms a separate set of slices for every object, pairs each slice with its assignment-weighted canonical position, and passes relative-position descriptors through a shared edge multilayer perceptron (MLP). A nodewise decoder maps the final hidden state to two real channels, which form one complex scattered-pressure value per canonical node.

Shared local parameters and the sum over source-object messages make the model equivariant to object permutations: consistently reordering the input objects reorders their predicted traces in the same way.

\subsection{Source-wise slice interactions}

The interaction module has eight heads, eight slices per head and parameters shared across its three applications. It uses assignments separate from those in the local Physics-Attention blocks. Let \(F_{ir}\in\R^C\) be the current feature at node \(r\) of object \(i\). For each head, an assignment \(\omega_{irs}^{(h)}\) is normalized over slices \(s\), and projected node values have width \(d_h=32\). The weighted token and its assignment-weighted position are
\begin{equation}
 z_{is}^{(h)}=\frac{\sum_r\omega_{irs}^{(h)}v_{ir}^{(h)}}{\mu_{is}^{(h)}},\qquad
 p_{is}^{(h)}=c_i+\frac{\sum_r\omega_{irs}^{(h)}\xi_r}{\mu_{is}^{(h)}},\qquad
 \mu_{is}^{(h)}=\sum_r\omega_{irs}^{(h)},
 \label{eq:sliceposition}
\end{equation}
where \(\xi_r\) is the local canonical coordinate and \(\mu_{is}^{(h)}\) is floored at a small positive value. Each token is paired with the assignment-weighted centroid of its canonical nodes, which defines the relative descriptor for a source--target slice pair.

For target slice \(s\) of object \(i\) and source slice \(t\) of object \(j\), write \(\delta_{ijst}^{(h)}=p_{is}^{(h)}-p_{jt}^{(h)}\in\R^3\) and \(r_{ijst}^{(h)}=|\delta_{ijst}^{(h)}|\). The implementation floors the distance before forming its unit direction. The nine-component relative descriptor, with distances normalized by \(R_c\), is
\begin{equation}
 e_{ijst}^{(h)}=\left[\frac{\delta_{ijst}^{(h)}}{R_c},\frac{r_{ijst}^{(h)}}{R_c},
 \frac{\delta_{ijst}^{(h)}}{r_{ijst}^{(h)}},\sin(kr_{ijst}^{(h)}),\cos(kr_{ijst}^{(h)})\right].
 \label{eq:edge}
\end{equation}
A shared \(9\to32\to2\) edge MLP is applied to each head's descriptor. Its parameters are shared across heads, while its outputs depend on the head-specific slice positions:
\begin{equation}
 (\beta_{ijst}^{(h)},\gamma_{ijst}^{(h)})=f_e(e_{ijst}^{(h)}),\qquad
 g_{ijst}^{(h)}=1+\tanh\gamma_{ijst}^{(h)}.
 \label{eq:edgeoutputs}
\end{equation}
The query, key and value maps \(Q,K,V\) are shared across heads and act on their respective head-specific tokens. The source-wise attention weights are
\begin{align}
 \ell_{ijst}^{(h)}&=\frac{(Qz_{is}^{(h)})^\top(Kz_{jt}^{(h)})}{\sqrt{d_h}}+\beta_{ijst}^{(h)},\label{eq:logit}\\
 A_{ijst}^{(h)}&=\operatorname{softmax}_{t}\bigl(\ell_{ijst}^{(h)}\bigr).\label{eq:sourceattention}
\end{align}
For fixed target object \(i\), target slice \(s\), and source object \(j\ne i\), the softmax is normalized over that source object's slices \(t\). Thus each source object contributes a separately normalized slice message; the target then sums these contributions:
\begin{equation}
 m_{is}^{(h)}=\sum_{j\ne i}\sum_t A_{ijst}^{(h)}g_{ijst}^{(h)}Vz_{jt}^{(h)}.
 \label{eq:sum}
\end{equation}
The target assignments return the messages to nodes. Head messages are concatenated, projected back to width \(C\) and added through the learned residual scale \(\alpha\):
\begin{equation}
 F_{ir}\leftarrow F_{ir}+\alpha L_{\mathrm{out}}
 \operatorname{concat}_{h}\left(\sum_s\omega_{irs}^{(h)}m_{is}^{(h)}\right).
 \label{eq:residual}
\end{equation}

For fixed hidden width and slice count, local processing scales linearly with the number of objects, whereas pairwise interactions scale quadratically.

\subsection{Prediction and operator-weighted supervision}
\label{sec:supervision}

The decoder predicts the complex scattered trace \(\widehat u_i^s\), from which the total trace follows as
\begin{equation}
 \widehat u_i^{\mathrm{tot}}=u_i^\inc+\widehat u_i^s.
 \label{eq:target}
\end{equation}
Training uses two scene-normalized losses. The weighted norm is \(\norm{v}_W^2=v^*Wv\), with the canonical-sphere area weights defined in Eq.~\eqref{eq:projection}:
\begin{align}
 \mathcal L_s&=\frac{\sum_i\norm{\widehat u_i^s-u_i^{s,\GT}}_W^2}
 {\sum_i\norm{u_i^{s,\GT}}_W^2+\varepsilon_s},\label{eq:scatloss}\\
 \mathcal L_q&=\frac{\sum_i\norm{Z_iP_Q(\widehat u_i^{\mathrm{tot}}-u_i^{\mathrm{tot},\GT})}_W^2}
 {\sum_i\norm{Z_iP_Qu_i^{\mathrm{tot},\GT}}_W^2+\varepsilon_q},\qquad
 \mathcal L=\frac{1}{B}\sum_{\mathrm{scene}}(\mathcal L_s+\lambda_q\mathcal L_q),
 \label{eq:qloss}
\end{align}
where \(B\) is the batch size and the stabilizers control small denominators. The auxiliary term is an operator-weighted flux error: for \(\delta u_i=\widehat u_i^{\mathrm{tot}}-u_i^{\mathrm{tot},\GT}\),
\(\norm{Z_iP_Q\delta u_i}_W^2=\delta u_i^*P_Q^*Z_i^*WZ_iP_Q\delta u_i\).
The auxiliary loss weights trace errors through the projected local flux operator \(Z_iP_Q\), while the scattered-field targets supervise scene coupling. Its reference flux is derived from the numerical reference trace and the same local response. Each scene is normalized before batch averaging, and training-set scales remain fixed at evaluation.

\section{Exterior-field recovery from canonical equivalent sources}
\label{sec:recovery}

\subsection{Pressure and flux as equivalent surface sources}

Adding the known incident trace to the predicted scattered pressure gives the total canonical trace \(\widehat u_i^{\mathrm{tot}}\). Applying the local response to this trace provides the flux representation for equivalent-source recovery. The continuous pressure--flux pair defines the surface identity below, and the numerical reconstruction uses its retained modal projection.

For an observation point outside all canonical balls, the outgoing Green function \(G_k(x,y)=e^{ik|x-y|}/(4\pi|x-y|)\) gives the canonical-surface representation
\begin{equation}
 u^s(x)=\sum_{i=1}^{N}\int_{\Gamma_i^c}
 \left[
 \partial_{n_y}G_k(x,y)\,u_i^{\mathrm{tot}}(y)
 -G_k(x,y)\,q_i^{\mathrm{tot}}(y)
 \right]\,dS_y.
 \label{eq:equivalentsources}
\end{equation}
The normal \(n_y\) points outward from each ball, giving the double-layer minus single-layer convention in Eq.~\eqref{eq:equivalentsources}. The incident wave and fields generated by other objects are regular inside each canonical ball; their complete pressure--flux pairs contribute zero to this exterior representation. The remaining contribution is the outgoing field generated by the enclosed object \cite{martin2006}. Summing over all objects gives the scene's scattered field.

\subsection{Recovering the outgoing coefficients}

The implementation projects the predicted total trace onto the retained \(Q\)-mode space:
\begin{equation}
 \alpha_i=P_Q\widehat u_i^{\mathrm{tot}},\qquad
 \eta_i=B_i\alpha_i,\qquad
 g_i^Q=\Phi_Q\alpha_i,\qquad q_i^Q=\Phi_Q\eta_i.
 \label{eq:recoveryprojection}
\end{equation}
Here \(g_i^Q\) and \(q_i^Q\) are the modal pressure and flux traces used for recovery. In contrast, the auxiliary training loss applies the full sampled response \(Z_iP_Q\) to the trace. The retained coefficients satisfy \(\alpha_i=Ja_i^{\mathrm{exc}}+Hb_i\) and \(\eta_i=J'a_i^{\mathrm{exc}}+H'b_i\). Multiplying the second relation by \(J\) and subtracting \(J'\) times the first eliminates the regular incoming contribution:
\begin{equation}
 J\eta_i-J'\alpha_i=(JH'-J'H)b_i.
 \label{eq:wronskianelimination}
\end{equation}
The diagonal radial Wronskian is~\cite[Section~10.50]{dlmf}
\(D_{\mathrm{Wr}}=JH'-J'H=i/(kR_c^2)I\),
where the radial derivatives include the factor \(k\). The outgoing coefficients are
\begin{equation}
 \widehat b_{i,lm}=
 \frac{j_l(kR_c)\eta_{i,lm}-j'_{l,r}\alpha_{i,lm}}
 {j_l(kR_c)h'_{l,r}-j'_{l,r}h_l^{(1)}(kR_c)},
 \label{eq:recoverb}
\end{equation}
where \(j'_{l,r}=kj_l'(kR_c)\) and \(h'_{l,r}=k[h_l^{(1)}]'(kR_c)\).

\subsection{Near-field and far-field evaluation}

At a point outside all canonical balls, define \(r_i=|x-c_i|\) and \(\widehat x_i=(x-c_i)/r_i\). The recovered near field is
\begin{equation}
 \widehat u^s(x)=\sum_{i=1}^{N}\sum_{l=0}^{L}\sum_{m=-l}^{l}
 \widehat b_{i,lm}h_l^{(1)}(kr_i)Y_l^m(\widehat x_i).
 \label{eq:near}
\end{equation}
All object contributions are evaluated at the same observation point and summed before comparison with a reference field.

For \(x=r\omega\) with \(|\omega|=1\), the asymptotic forms
\(h_l^{(1)}(kr_i)\sim(-i)^{l+1}e^{ikr_i}/(kr_i)\)~\cite[Eq.~10.52.4]{dlmf} and
\(r_i=r-\omega\cdot c_i+O(r^{-1})\) give
\(u^s(r\omega)=e^{ikr}f(\omega)/r+O(r^{-2})\).
The far-field amplitude \(f=u^\infty\) is therefore reconstructed as
\begin{equation}
 \widehat f(\omega)=\frac1k\sum_{i=1}^{N}e^{-ik\omega\cdot c_i}
 \sum_{l=0}^{L}\sum_{m=-l}^{l}(-i)^{l+1}Y_l^m(\omega)\widehat b_{i,lm}.
 \label{eq:far}
\end{equation}
The translation phase preserves interference between radiation from different centres. For fixed local operators and object positions, this reconstruction is a linear trace-to-field map.

For a fixed scene, let \(P_N=\operatorname{diag}(P_Q,\ldots,P_Q)\) project the concatenated canonical traces, and let \(\mathcal E_Q^X\) evaluate and sum all objects' outgoing expansions on an exterior observation set \(X\). The retained trace-to-field operator is
\begin{equation}
 \mathcal R_Q^X=\mathcal E_Q^X\operatorname{diag}_{i=1}^N
 \!\left[D_{\mathrm{Wr}}^{-1}(JB_i-J')\right]P_N.
 \label{eq:recoveryoperator}
\end{equation}
With canonical-interface and observation-set quadrature weights \(W_N\) and \(W_X\), linearity gives
\begin{equation}
 \norm{\widehat u_X^s-u_{\mathrm{ref},X}^s}_{W_X}
 \leq \norm{\mathcal R_Q^X}_{W_N\to W_X}
 \norm{\widehat u^{\mathrm{tot}}-u^{\mathrm{tot,ref}}}_{W_N}
 +E_{\mathrm{repr}}(X),
 \label{eq:fielderrordecomp}
\end{equation}
where \(E_{\mathrm{repr}}(X)=\norm{\mathcal R_Q^Xu^{\mathrm{tot,ref}}-u_{\mathrm{ref},X}^s}_{W_X}\). The first term bounds propagation of the canonical-trace prediction error, and the second is the retained-mode representation discrepancy relative to the numerical reference. Dividing both sides by the same reference-field norm gives the corresponding relative-error bound. The reported near- and far-field metrics use their respective planar-grid and spherical-quadrature norms. Appendix~\ref{app:theory} gives the detailed derivation and distinguishes trace-prediction error from local-operator perturbations and finite-mode representation error.

\section{Experiments and results}
\label{sec:results}

\subsection{Experimental setup}

Numerical experiments are conducted on three-dimensional acoustic scattering by sound-hard obstacles at a fixed wavenumber of \(k=2\). Each scene contains two to four scatterers.

The training data comprise two datasets. The first contains 10,000 scenes constructed from spheres, ellipsoids, and boxes. The second contains 10,000 scenes generated from a library of 11 CAD geometries, including unperforated, single-bore, and double-bore ellipsoids with rounded bore edges. Both datasets include variations in object arrangements, orientations, and incident excitations. An independent validation set of 1,000 scenes is provided for each dataset.

The \textbf{Same-geometry structure} test set evaluates generalization to held-out instances of the sphere, ellipsoid, and box families represented during training. It contains 1,000 scenes constructed from base geometries that are excluded from the training library.

\begin{figure}[!htbp]
\centering
\includegraphics[width=\textwidth]{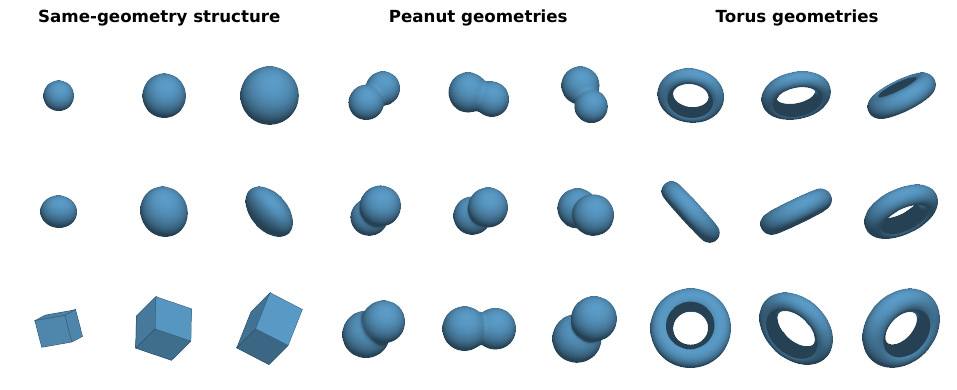}
\caption{Geometries from the evaluated test sets: Same-geometry structure, Peanut geometries and Torus geometries.}
\label{fig:testgeometries}
\end{figure}

Two additional test sets assess generalization to unseen geometry families. The \textbf{Peanut geometries} set contains 500 scenes assembled from 24 nonconvex Peanut geometries, while the \textbf{Torus geometries} set contains 500 scenes assembled from 24 Torus geometries. Each test set is divided into 250 mixed scenes, each containing one scatterer from the corresponding unseen family, and 250 scenes composed entirely of scatterers from that family. Figure~\ref{fig:testgeometries} illustrates representative base geometries. Different scenes are generated by combining and rotating these geometric instances.

The canonical enclosing spheres have a radius of \(R_c=1.5\) and are discretized using 128 Fibonacci sampling points. ESPC-Transolver takes 64 probe responses corresponding to spherical-wave modes truncated at \(l=7\), whereas the reference-field calculations use 100 input modes truncated at \(l=9\). The incident fields are unit-amplitude plane waves with propagation directions sampled uniformly over the unit sphere. Further details of the numerical discretization and model configurations are provided in Appendix~\ref{app:methodconfig}.

\subsection{Experiment 1: Comparison of response representations under a fixed compositional backbone}

This experiment evaluates the effect of the conditioning representation while keeping the ESPC-Transolver backbone unchanged.

\textbf{Transolver-Geo-matched} replaces the canonical scattering-response channels with radial-occupancy features evaluated along the canonical sampling directions. \textbf{Transolver-Geo-strong} uses surface points and their associated normals to construct both local neighborhood features and a pooled global representation of each scatterer.

Both geometry-conditioned models provide 128 conditioning channels. Unlike ESPC-Transolver, which uses canonical scattering responses as input, these models directly condition their predictions on geometric features. All three compositional models use the same operator-weighted training loss and the same local response operators for exterior-field reconstruction. The encoder architectures and parameter counts are detailed in Appendix~\ref{app:methodconfig}.

\begin{figure}[!htbp]
\centering
\includegraphics[width=\textwidth]{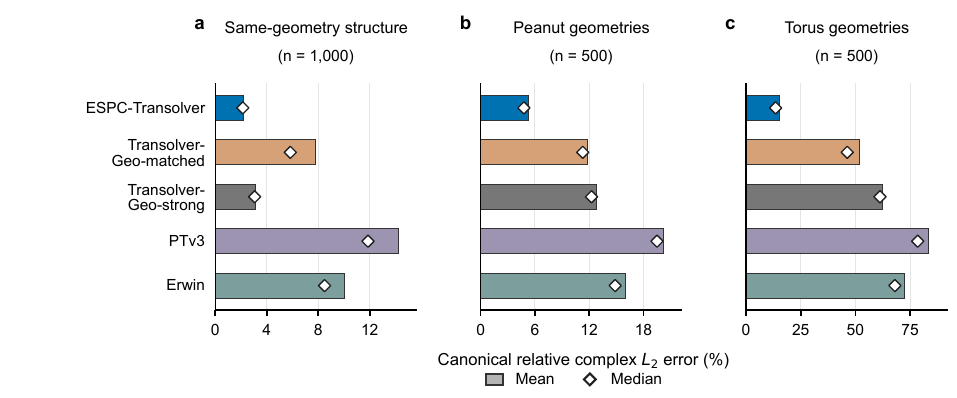}
\caption{Canonical scattered-field errors for all methods in Table~\ref{tab:representation20k}. Bars show scene means and diamonds show medians. The three panels use separate horizontal scales.}
\label{fig:representation20k}
\end{figure}

The prediction results (Table~\ref{tab:representation20k} and Fig.~\ref{fig:representation20k}) show that ESPC-Transolver achieves lower errors than all other models on the test set comprising held-out instances of familiar geometry families. Among the geometry-conditioned models, Transolver-Geo-strong outperforms Transolver-Geo-matched on familiar families, with mean errors of 3.10\% and 7.77\%, respectively, but their ordering reverses on Peanut and Torus. The advantage of richer surface encoding within familiar families therefore does not extend to these unseen families.

Canonical scattering-response conditioning yields lower mean and median errors than both geometry-conditioned alternatives on the complete Peanut and Torus test sets. This pattern is consistent with the shared excitation basis providing physically comparable response channels across geometry families.

\begin{table}[!htbp]
\centering
\small
\renewcommand{\arraystretch}{1.02}
\setlength{\tabcolsep}{5pt}
\caption{Canonical scattered-field relative complex \(L_2\) error (\%) for ESPC-Transolver and four geometry-conditioned predictors.}
\label{tab:representation20k}
\begin{tabular}{@{}>{\raggedright\arraybackslash}p{3.1cm}>{\raggedright\arraybackslash}p{5.2cm}rr@{}}
\toprule
Test set & Method & Mean & Median \\
\midrule
Same-geometry structure & \textbf{ESPC-Transolver} & \textbf{2.21} & \textbf{2.14} \\
 & Transolver-Geo-matched & 7.77 & 5.83 \\
 & Transolver-Geo-strong & 3.10 & 3.06 \\

 & PTv3 & 14.20 & 11.87 \\
 & Erwin & 10.05 & 8.49 \\
\addlinespace[5pt]
Peanut geometries & \textbf{ESPC-Transolver} & \textbf{5.27} & \textbf{4.78} \\
 & Transolver-Geo-matched & 11.76 & 11.29 \\
 & Transolver-Geo-strong & 12.80 & 12.26 \\

 & PTv3 & 20.21 & 19.50 \\
 & Erwin & 16.00 & 14.87 \\
\addlinespace[5pt]
Torus geometries & \textbf{ESPC-Transolver} & \textbf{15.44} & \textbf{13.49} \\
 & Transolver-Geo-matched & 52.03 & 46.27 \\
 & Transolver-Geo-strong & 62.29 & 61.23 \\

 & PTv3 & 83.56 & 78.47 \\
 & Erwin & 72.52 & 68.04 \\
\bottomrule
\end{tabular}
\end{table}

\begin{figure}[p]
\centering
\includegraphics[width=\textwidth]{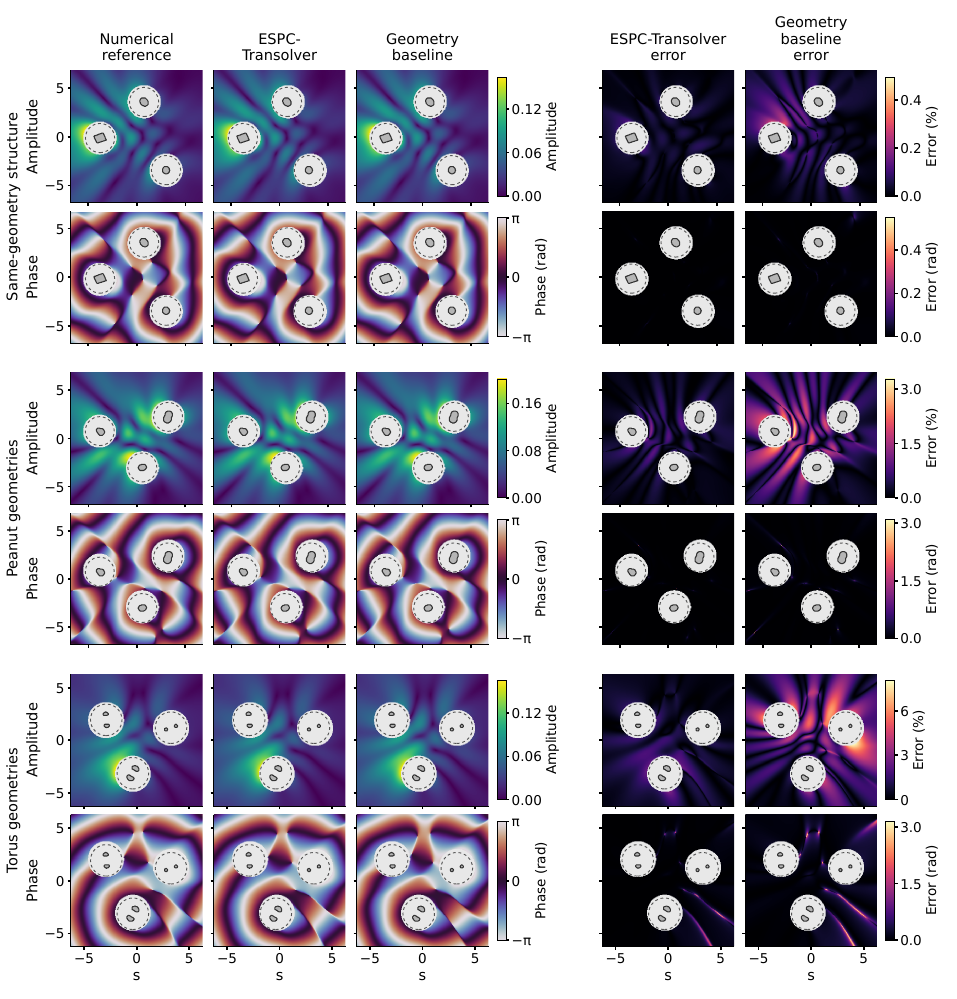}
\caption{Near-field amplitude and phase for the three test scenes.}
\label{fig:representationfield20k}
\par\vspace{4pt}
\includegraphics[width=\textwidth]{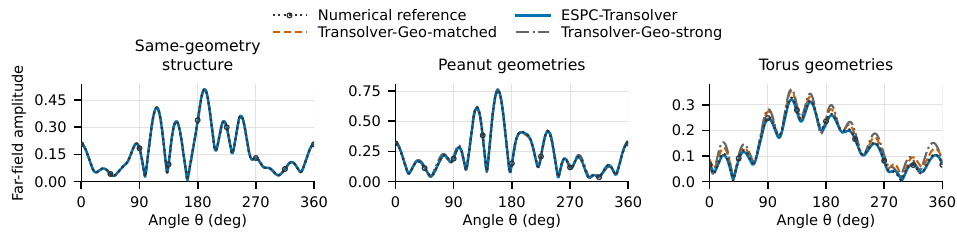}
\caption{Far-field amplitudes along a great circle for the numerical reference, ESPC-Transolver, Transolver-Geo-matched and Transolver-Geo-strong in the same scenes.}
\label{fig:representationfar20k}
\end{figure}

Figure~\ref{fig:representationfield20k} compares the exterior-field reconstructions obtained with ESPC-Transolver and the geometry-conditioned model yielding the lower mean canonical scattered-field error on each test set. Transolver-Geo-strong is selected as the geometric baseline for Same-geometry structure, while Transolver-Geo-matched is used for Peanut and Torus.

For each test set, a representative scene containing three scatterers is selected to visualize the amplitude and phase distributions of the reconstructed scattered field. Each pair of rows presents the corresponding amplitude and phase results for one scene.

\clearpage

\begin{table}[H]
\centering
\small
\renewcommand{\arraystretch}{1.02}
\setlength{\tabcolsep}{3pt}
\caption{Recovery errors for the three illustrated scenes. Near-field complex and amplitude errors and far-field complex errors are relative \(L_2\) values (\%); phase mean absolute error (MAE) is in radians.}
\label{tab:representationfields20k}
\begin{tabular}{@{}>{\raggedright\arraybackslash}p{2.8cm}>{\raggedright\arraybackslash}p{4.6cm}rrrr@{}}
\toprule
Test geometry & Method & \shortstack{Near\\\(E_s\)} & \shortstack{Near\\\(E_{\rm amp}\)} & \shortstack{Phase\\MAE} & \shortstack{Far\\\(E_s\)} \\
\midrule
Same-geometry structure & \textbf{ESPC-Transolver} & \textbf{0.15} & \textbf{0.11} & \textbf{0.0015} & \textbf{0.13} \\
 & Transolver-Geo-matched & 0.74 & 0.43 & 0.0091 & 0.73 \\
 & Transolver-Geo-strong & 0.32 & 0.23 & 0.0036 & 0.29 \\
 & PTv3 & 2.22 & 1.61 & 0.0237 & 2.10 \\
 & Erwin & 1.37 & 0.77 & 0.0190 & 1.26 \\
\addlinespace[5pt]
Peanut geometries & \textbf{ESPC-Transolver} & \textbf{1.20} & \textbf{0.83} & \textbf{0.0118} & \textbf{1.22} \\
 & Transolver-Geo-matched & 3.22 & 2.36 & 0.0277 & 3.32 \\
 & Transolver-Geo-strong & 5.77 & 3.79 & 0.0661 & 6.12 \\
 & PTv3 & 5.47 & 3.44 & 0.0567 & 5.80 \\
 & Erwin & 4.26 & 2.68 & 0.0453 & 4.48 \\
\addlinespace[5pt]
Torus geometries & \textbf{ESPC-Transolver} & \textbf{4.35} & \textbf{3.23} & \textbf{0.0713} & \textbf{4.42} \\
 & Transolver-Geo-matched & 12.80 & 8.78 & 0.1772 & 11.85 \\
 & Transolver-Geo-strong & 15.36 & 12.94 & 0.1416 & 15.68 \\
 & PTv3 & 26.31 & 22.75 & 0.2206 & 26.53 \\
 & Erwin & 13.79 & 10.69 & 0.1836 & 13.08 \\
\bottomrule
\end{tabular}
\end{table}

ESPC-Transolver achieves lower near-field complex-valued reconstruction errors than both geometry-conditioned models in all three representative scenes (Table~\ref{tab:representationfields20k}). The far-field reconstruction errors are also lower for ESPC-Transolver; Fig.~\ref{fig:representationfar20k} compares the angular patterns of all three predictors with the numerical reference.

These results indicate that canonical scattering-response conditioning improves not only the accuracy of interface-field predictions but also the reconstruction of exterior scattered fields, including the interference patterns arising from multiple wave interactions.

\subsection{Experiment 2: Comparison of neural network architectures}

ESPC-Transolver models inter-object coupling by alternating object-local feature processing with relative-position-aware interactions. In contrast, Transolver employs Physics-Attention. The PIBNet multiscale graph backbone is adapted to predict canonical-interface scattered pressure from the response channels defined above~\cite{marsal2025}. All five models in this architectural comparison use the same canonical scattering-response inputs. PTv3~\cite{wuptv3} and Erwin~\cite{zhdanov2025} retain their hierarchical point-processing backbones and use pressure-regression heads specified in Appendix~\ref{app:methodconfig}.

For the quantitative comparison, Transolver serves as the baseline, whereas ESPC-Transolver consists of five object-local processing blocks with three interleaved applications of a shared interaction module. Detailed architectural configurations are provided in Appendix~\ref{app:methodconfig}.

As reported in Table~\ref{tab:architecture20k}, ESPC-Transolver achieves the lowest mean and median prediction errors among all evaluated architectures. Compared with Transolver, ESPC-Transolver reduces the mean canonical scattered-field error from 4.97\% to 2.21\% on the Same-geometry structure test set and from 7.46\% to 5.27\% on the Peanut test set. The difference is smaller on Torus, where the corresponding errors are 16.05\% and 15.44\%, respectively. The advantage of the coupling-aware architecture over Transolver is therefore concentrated in the familiar families and Peanut, whereas the two response-conditioned Transformers are closer on Torus. Together with Table~\ref{tab:representation20k}, this result distinguishes the benefit of canonical-response conditioning on Torus from the smaller improvement associated with the coupling-aware architecture.

\begin{figure}[p]
\centering
\includegraphics[width=\textwidth]{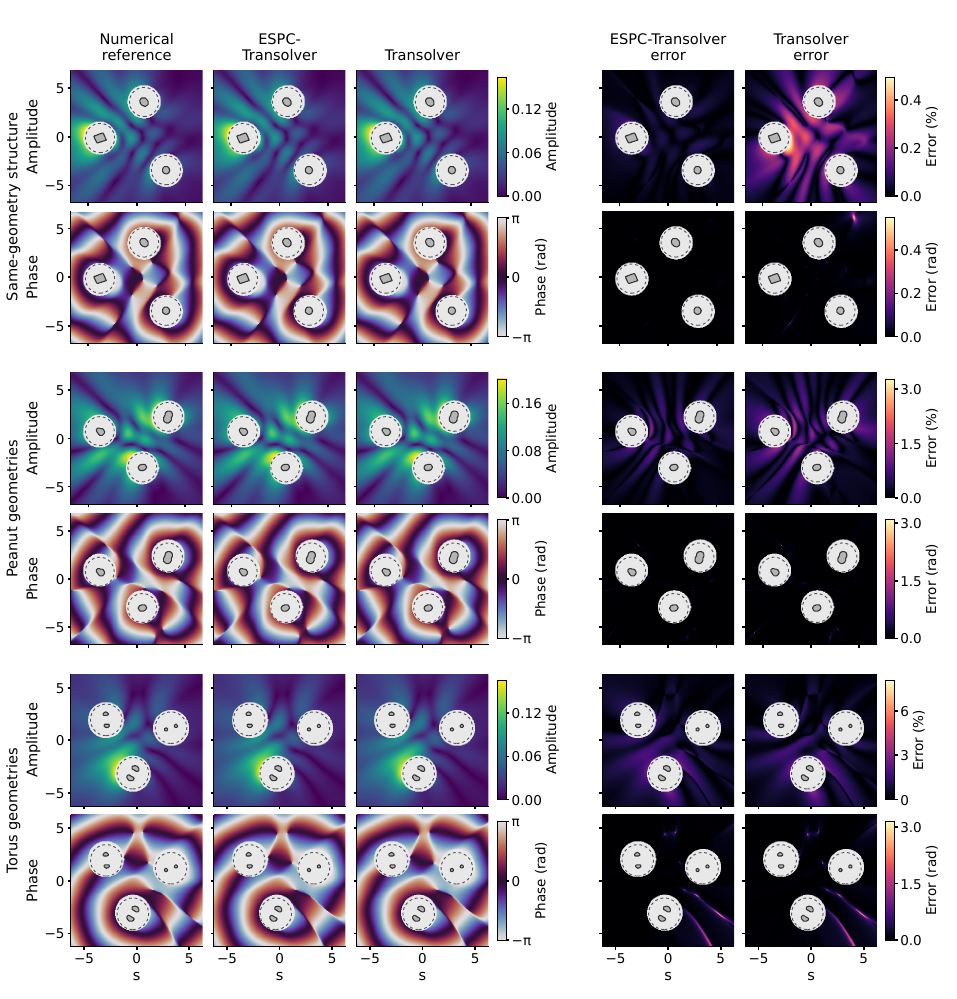}
\caption{Near-field amplitude and phase for the three test scenes.}
\label{fig:architecturefield20k}
\par\vspace{4pt}
\includegraphics[width=\textwidth]{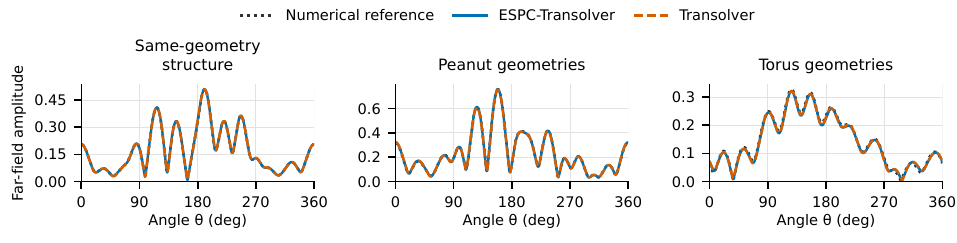}
\caption{Far-field amplitudes along a great circle for ESPC-Transolver and Transolver in the same three scenes.}
\label{fig:architecturefar20k}
\end{figure}

\begin{table}[!htbp]
\centering
\small
\renewcommand{\arraystretch}{1.02}
\setlength{\tabcolsep}{5pt}
\caption{Interaction architectures with canonical scattering-response inputs. Values are canonical scattered-field relative complex \(L_2\) errors (\%).}
\label{tab:architecture20k}
\begin{tabular}{@{}>{\raggedright\arraybackslash}p{3.1cm}>{\raggedright\arraybackslash}p{5.2cm}rr@{}}
\toprule
Test set & Method & Mean & Median \\
\midrule
Same-geometry structure & \textbf{ESPC-Transolver} & \textbf{2.21} & \textbf{2.14} \\
 & Transolver & 4.97 & 4.11 \\
 & PTv3 & 14.33 & 11.91 \\
 & Erwin & 15.44 & 15.53 \\
 & PIBNet & 13.43 & 13.92 \\
\addlinespace[5pt]
Peanut geometries & \textbf{ESPC-Transolver} & \textbf{5.27} & \textbf{4.78} \\
 & Transolver & 7.46 & 7.12 \\
 & PTv3 & 17.91 & 16.68 \\
 & Erwin & 17.49 & 17.17 \\
 & PIBNet & 16.79 & 17.13 \\
\addlinespace[5pt]
Torus geometries & \textbf{ESPC-Transolver} & \textbf{15.44} & \textbf{13.49} \\
 & Transolver & 16.05 & 15.18 \\
 & PTv3 & 29.03 & 27.83 \\
 & Erwin & 26.18 & 25.76 \\
 & PIBNet & 25.76 & 24.51 \\
\bottomrule
\end{tabular}
\end{table}
Figures~\ref{fig:architecturefield20k} and \ref{fig:architecturefar20k} show field illustrations for ESPC-Transolver and Transolver in the same three scenes used in Experiment~1. Their amplitude and phase maps and far-field curves illustrate the recovered pressure patterns. Table~\ref{tab:architecture20k} provides the quantitative comparison with Transolver.

Table~\ref{tab:inference100} reports the computational time required to predict 100 scattering problems. For neural network models, the measurements include only the forward inference time.

ESPC-Transolver requires 1.088~s for 100 predictions, compared with 1.264~s for Transolver-Geo-matched and 1.239~s for Transolver-Geo-strong. Transolver completes the same number of predictions in 0.540~s.

\begin{table}[!htbp]
\centering
\small
\renewcommand{\arraystretch}{1.02}
\setlength{\tabcolsep}{5pt}
\caption{Time for 100 three-object scenes: online numerical solution and neural forward inference, excluding local-response precomputation and training.}
\label{tab:inference100}
\begin{tabular}{@{}>{\raggedright\arraybackslash}p{5.7cm}>{\raggedright\arraybackslash}p{6.0cm}r@{}}
\toprule
Method & Input representation & Time (s) \\
\midrule
T-matrix-based Foldy--Lax solver & --- & 33.089 \\
\addlinespace[3pt]
ESPC-Transolver & Canonical scattering responses & 1.088 \\
Transolver-Geo-matched & Radial occupancy & 1.264 \\
Transolver-Geo-strong & Surface points and normals & 1.239 \\
PTv3 & Canonical scattering responses & 3.430 \\
Erwin & Canonical scattering responses & 2.089 \\
Transolver & Canonical scattering responses & \textbf{0.540} \\
PTv3 & Surface points and normals & 3.665 \\
Erwin & Surface points and normals & 1.911 \\
\bottomrule
\end{tabular}
\end{table}

Together, Tables~\ref{tab:architecture20k} and~\ref{tab:inference100} illustrate the trade-off between prediction accuracy and inference cost for Transolver and coupling-aware object-local processing. ESPC-Transolver achieves lower scattered-field prediction errors through object-local feature processing and relative-position-aware interactions, at the expense of increased inference time relative to Transolver.

\section{Conclusion}
This work introduced ESPC-Transolver, which separates object-specific scattering responses from scene-dependent interactions to predict acoustic multiple scattering. Canonical scattering responses provide a physically consistent representation across geometries, while shared object-local blocks and relative-position-aware attention model inter-object coupling. Across three-dimensional sound-hard scattering problems, ESPC-Transolver achieves lower interface and canonical scattered-field prediction errors than the evaluated geometry-conditioned baselines and alternative network architectures, including on held-out examples of familiar families and previously unseen Peanut and Torus geometries. Near- and far-field reconstructions further show that the predicted interface pressure and normal flux recover exterior pressure distributions and multiple-scattering interference patterns, demonstrating the effectiveness of this approach for cross-geometry scattering prediction.

\section*{Acknowledgments} 
This work is supported by NSFC (Grant No. 12571399 and 12101511) and Fundamental Research Funds for the Central Universities (Grant No. JBK22YJ37).

\section*{Declaration of generative AI and AI-assisted technologies in the manuscript preparation process}
During the preparation of this work, the authors used ChatGPT by OpenAI to support language editing, improve readability, and obtain feedback on manuscript structure and phrasing. The authors reviewed, edited, and verified the content as needed and take full responsibility for the content of the article.

\appendix
\begingroup
\small
\setlength{\parskip}{2pt}
\setlength{\abovedisplayskip}{6pt plus 2pt minus 2pt}
\setlength{\belowdisplayskip}{6pt plus 2pt minus 2pt}
\setlength{\abovedisplayshortskip}{3pt plus 1pt minus 1pt}
\setlength{\belowdisplayshortskip}{3pt plus 1pt minus 1pt}
\section{Network configurations and numerical implementation}
\label{app:methodconfig}
\label{app:numerics}

Table~\ref{tab:modelconfig} summarizes the network configurations and total trainable parameter counts of all evaluated models, including input encoders and output prediction heads. Parameters in shared modules are counted only once.

\begin{table}[!htbp]
\centering
\footnotesize
\renewcommand{\arraystretch}{1.02}
\setlength{\tabcolsep}{4pt}
\caption{Network configurations and total trainable parameters. Layer widths, depths, and attention head counts are listed in stage order.}
\label{tab:modelconfig}
\begin{tabular}{@{}p{4.0cm}p{8.7cm}r@{}}
\toprule
Method & Configuration & Parameters \\
\midrule
ESPC-Transolver & Five object-local blocks with a hidden width of 256, 8 attention heads, 32 local slices per head, and an MLP ratio of 1. A shared interaction module with 8 heads and 8 slices per head is applied after each of the first three local blocks. & 2.077M \\
\addlinespace[3pt]
Transolver-Geo-matched & ESPC-Transolver backbone with 128 radial-occupancy channels replacing the canonical scattering-response channels. & 2.092M \\
\addlinespace[3pt]
Transolver-Geo-strong & ESPC-Transolver backbone augmented with local and global surface-point encoders. & 2.144M \\
\addlinespace[3pt]
Transolver & Physics-Attention with a hidden width of 256, 8 attention heads, 32 slices per head, an MLP ratio of 1, and zero dropout. & 1.541M \\
\addlinespace[3pt]
PIBNet & Three-level canonical-node graph with a base hidden width of 64 and a width scaling factor of 2 between levels. The first and last top-level processors use 3 steps each, while the bottom-level processor uses 6 steps. SiLU activation and an edge feature dimension of 134 are used. & 5.221M \\
\addlinespace[3pt]
PTv3 & Encoder depths \((2,2,2,2)\), widths \((32,64,128,256)\), and attention heads \((2,4,8,16)\); decoder depths \((2,2,2)\), widths \((64,64,128)\), and attention heads \((4,4,8)\). Strides \((2,2,2)\), patch size 128, and grid size 0.16. Z, Z-trans, Hilbert, and Hilbert-trans orders are used. Drop-path is set to zero; order shuffling, Flash Attention, and relative positional encoding are disabled. & 9.699M \\
\addlinespace[3pt]
Erwin & Encoder widths \((32,64,128,256)\), attention heads \((2,4,8,16)\), and depths \((2,2,2,2)\); decoder attention heads \((2,4,8)\) and depths \((2,2,2)\). Strides \((2,2,2)\), ball sizes \((32,32,32,32)\), rotation setting 45, MLP ratio 4, and zero message-passing steps. & 3.753M \\
\bottomrule
\end{tabular}
\end{table}

\subsection{Geometric feature encoding and network implementation}

Radial-occupancy features are evaluated at the radial sampling locations
\[
 \rho_s=\frac{0.9(s-\frac{1}{2})}{128},\qquad s=1,\ldots,128.
\]
Feature normalization uses the means and variances computed from the training data, with a variance floor of \(10^{-4}\).

The surface encoders operate on 2,048 area-sampled surface points and their associated normals. Local and global surface features are extracted using two independent multilayer perceptrons (MLPs), each with dimensions \(6\rightarrow64\rightarrow128\).
The local encoder processes the 16 nearest neighbors and aggregates their features through max pooling. The global encoder processes the entire point cloud and combines max and mean pooling to construct a whole-object representation.

The resulting 128-dimensional local features and 256-dimensional global features are concatenated and projected onto a 128-dimensional conditioning space through a \(384\rightarrow128\) linear projection.

The regression heads of PTv3 and Erwin use the architectures \(64\rightarrow64\rightarrow2\) and \(32\rightarrow64\rightarrow2\), respectively.

PIBNet constructs its graph using triangulation edges, a depth-two binary tree, and reciprocal connections between downsampling and upsampling levels. Distant edges are generated using an edge ratio of 0.1, two candidates, and a sampling period of 6. Self-loops and duplicate edges are removed.

\subsection{Training configuration}

All models are trained using the AdamW optimizer for 200 epochs with a batch size of 16. The learning rate follows a cosine annealing schedule, decreasing from \(10^{-4}\) to \(10^{-7}\). Weight decay is set to \(10^{-5}\), and gradients are clipped to a maximum norm of 1.

The auxiliary loss weight is set to 0.1, except for PTv3 and Erwin, for which it is zero.

The stabilization parameter in the scattered-field loss is defined as
\[
 \varepsilon_s=10^{-6}\sum_i\left\|u_i^{\mathrm{inc}}\right\|_W^2.
\]
The stabilization parameter \(\varepsilon_q\) in the flux loss is set to \(10^{-6}\) times the median flux energy over the training scenes.

Each response degree of freedom and the global scattered-field output are normalized using uncentered root-mean-square (RMS) scales computed from the training data. The original physical scales are restored during evaluation.

\subsection{Numerical discretization and implementation details}

The numerical discretization employs graded, curved quadratic finite-element (P2 FEM) meshes.

The retained modal dimensions are
\[
 Q_{\mathrm{in}}=100,\qquad Q_{\mathrm{out}}=121,
\]
with
\[
 E_Q\in\mathbb{R}^{100\times64},\qquad
 \Phi_{\mathrm{out}}\in\mathbb{C}^{128\times121},
\]
as defined in Eq.~\eqref{eq:fluxprojection}.

Spherical-harmonic modes are ordered by increasing degree \(l\), with the order \(m\) ranging from \(-l\) to \(l\) for each degree. The phase convention follows that of the numerical library used in the implementation.

\subsection{Error metrics}

The scene-wise relative scattered-field error is defined as
\[
 E_s=\frac{\left\|\widehat{\mathbf{u}}^{\,s}-\mathbf{u}_{\mathrm{ref}}^{\,s}\right\|_{W_N}}
 {\left\|\mathbf{u}_{\mathrm{ref}}^{\,s}\right\|_{W_N}},
\]
where \(\widehat{\mathbf{u}}^{\,s}\) and \(\mathbf{u}_{\mathrm{ref}}^{\,s}\) denote the predicted and reference scattered fields, respectively, and \(\|\cdot\|_{W_N}\) denotes the corresponding weighted norm.

The relative complex-field and amplitude errors for exterior-field reconstruction are normalized by the corresponding reference-field norms. Near-field errors are evaluated using uniform planar sampling weights, while far-field errors are computed using spherical quadrature weights.

The phase mean absolute error (MAE) is defined as
\[
 E_{\mathrm{phase}}=\operatorname{mean}_{x\in X_{\mathrm{valid}}}
 \left|\arg\left(\widehat{u}^{s}(x)\overline{u_{\mathrm{ref}}^{s}(x)}\right)\right|,
\]
where \(X_{\mathrm{valid}}\) denotes the set of valid exterior evaluation points. Phase errors are reported in radians.

\section{Supporting derivations}
\label{app:theory}

\subsection{Finite-dimensional projection and coupling consistency}

Set \(C_i^{(Q)}=C_iE_Q\) and partition its output modes at \(l=7\). Since \(\Phi_{\mathrm{out}}=[\Phi_Q,\Phi_{>7}]\),
\begin{equation}
 B_i=P_QZ_i=(C_i^{(Q)})_{\leq7}+P_Q\Phi_{>7}(C_i^{(Q)})_{>7}.
 \label{eq:discreteleakage}
\end{equation}
The second term records the finite-quadrature projection of higher output modes; the implemented conversion uses this sampled \(B_i\).

Let \(\mathsf T=\operatorname{diag}(T_i^\out)\) and let \(\mathsf U\) denote the off-diagonal translation operator. The retained solution satisfies \((I-\mathsf T\mathsf U)b_Q=\mathsf T a^\inc\). For a predicted total trace, set \(\alpha_i=P_Q\widehat u_i^{\mathrm{tot}}\), \(\eta_i=B_i\alpha_i\), and \(D_{\mathrm{Wr}}=JH'-J'H\). Eliminating regular and outgoing radial factors gives
\begin{equation}
 \widehat a_i^\exc=D_{\mathrm{Wr}}^{-1}(H'\alpha_i-H\eta_i),\qquad
 \widehat b_i=D_{\mathrm{Wr}}^{-1}(J\eta_i-J'\alpha_i).
 \label{eq:recoveredcoefficients}
\end{equation}
These coefficients satisfy \(\widehat b_i=T_i^\out\widehat a_i^\exc\). Stacking the objects and defining \(r^\mathrm{coup}=\widehat a^\exc-a^\inc-\mathsf U\widehat b\) yields
\begin{equation}
 r^\mathrm{FL}=(I-\mathsf T\mathsf U)\widehat b-\mathsf T a^\inc
 =\mathsf T r^\mathrm{coup},\qquad
 \norm{\widehat b-b_Q}\leq\norm{(I-\mathsf T\mathsf U)^{-1}}\norm{r^\mathrm{FL}}.
 \label{eq:couplingbound}
\end{equation}
For nonsingular \(I-\mathsf T\mathsf U\), the bound follows by subtracting the reference equation and relates recovered coefficients to physical coupling. Training uses the losses in Section~\ref{sec:supervision}.

\subsection{Local-response and coupled-system sensitivity}

For fixed centres and illumination, let tildes denote a perturbed local response and its corresponding retained solution, with both coupled systems nonsingular. Write \(A_i=H'-B_iH\). Subtracting the local conversion equations and the coupled equations, respectively, gives
\begin{align}
 T_i^\out-\widetilde T_i^\out
 &=A_i^{-1}(B_i-\widetilde B_i)(J+H\widetilde T_i^\out),\label{eq:localperturbation}\\
 (I-\mathsf T\mathsf U)(b_Q-\widetilde b_Q)
 &=(\mathsf T-\widetilde{\mathsf T})(a^\inc+\mathsf U\widetilde b_Q).
 \label{eq:identity}
\end{align}
Using \(B_i-\widetilde B_i=P_Q(Z_i-\widetilde Z_i)\) and submultiplicativity in compatible induced norms gives
\begin{align}
 \norm{T_i^\out-\widetilde T_i^\out}
 &\leq\norm{A_i^{-1}}\norm{P_Q}\norm{Z_i-\widetilde Z_i}
 \bigl(\norm{J}+\norm{H}\norm{\widetilde T_i^\out}\bigr),\label{eq:zperturbation}\\
 \norm{b_Q-\widetilde b_Q}
 &\leq\norm{(I-\mathsf T\mathsf U)^{-1}}
 \norm{\mathsf T-\widetilde{\mathsf T}}
 \bigl(\norm{a^\inc}+\norm{\mathsf U}\norm{\widetilde b_Q}\bigr).
 \label{eq:physicalbound}
\end{align}
The estimates separate local-response sensitivity from coupled amplification, whose spectral-norm factor is \(1/\sigma_{\min}(I-\mathsf T\mathsf U)\).

\subsection{Trace-prediction and exterior-representation errors}

At fixed local responses, let \(K_i=D_{\mathrm{Wr}}^{-1}(JB_i-J')P_Q\). Equation~\eqref{eq:recoverb} gives
\begin{equation}
 \widehat b_i-b_i^\mathrm{trace}=K_i
 (\widehat u_i^\mathrm{tot}-u_i^{\mathrm{tot},\GT}),
 \label{eq:tracebound}
\end{equation}
where \(b_i^\mathrm{trace}=K_i u_i^{\mathrm{tot},\GT}\) applies the same retained recovery to the reference trace. The exterior map is therefore \(\mathcal R_Q^X=\mathcal E_Q^X\operatorname{diag}(K_i)\), with \(\mathcal E_Q^X\) summing all objects' outgoing expansions. Adding and subtracting the recovered reference trace gives the exact decomposition
\begin{equation}
 \widehat u_X^s-u_{\mathrm{ref},X}^s
 =\mathcal R_Q^X(\widehat u^\mathrm{tot}-u^{\mathrm{tot},\GT})
 +\bigl(\mathcal R_Q^X u^{\mathrm{tot},\GT}-u_{\mathrm{ref},X}^s\bigr).
 \label{eq:appendixfielddecomposition}
\end{equation}
Taking weighted norms gives Eq.~\eqref{eq:fielderrordecomp}, separating neural trace error from finite-mode recovery and numerical-representation discrepancy. Equations~\eqref{eq:zperturbation}--\eqref{eq:physicalbound} describe local-response perturbations.

\endgroup
\let\manuscriptbibliography\thebibliography
\renewcommand{\thebibliography}[1]{\manuscriptbibliography{#1}\setlength{\baselineskip}{10pt}\setlength{\itemsep}{0pt}\setlength{\parsep}{0pt}}

\end{document}